\documentclass[a4paper,12pt,fleqn]{article}

\usepackage{amsfonts,latexsym,amsmath}

\newtheorem{theorem}{Theorem}
\newtheorem{proposition}[theorem]{Proposition}
\newtheorem{lemma}[theorem]{Lemma}

\newtheorem{definition}[theorem]{Definition}

\def\findem{\rule{0.2cm}{0.2cm}}

\def\bhyp#1{\begin{equation}\label{#1}\begin{array}{c}}	
\def\ehyp{\end{array}\end{equation}}

\def\bpro#1{\begin{equation}\label{#1}\left\{\begin{array}{l@{\qquad}l}}
 
\def\epro{\end{array}\right.\end{equation}}

\def\N{\mathbb{N}}
\def\Z{\mathbb{Z}}

\def\R{\mathbb{R}}

\def\<{\langle}
\def\>{\rangle}

\def\ds{\displaystyle} 
\def\div{{\rm div}}
\def\refe#1{(\ref{#1})}

\def\ocirc#1{\ifmmode\setbox0=\hbox{$#1$}\dimen0=\ht0
    \advance\dimen0 by1pt\rlap{\hbox to\wd0{\hss\raise\dimen0
    \hbox{\hskip.2em$\scriptscriptstyle\circ$}\hss}}#1\else
    {\accent"17 #1}\fi} 

\def\sgn{\mathrm{sgn}}
\def\eps{\varepsilon}
\def\qed{\findem}

\def\vg{\mathbf{v}_g}
\def\<{\langle}
\def\>{\rangle}

\def\S{\mathbb{S}}
\def\T{\mathbb{T}}

\begin{document}

\title{Long-time behavior in scalar conservation laws}
\author{A. Debussche and J. Vovelle}
\maketitle

\begin{abstract} We consider the long-time behavior of the entropy solution of a first-order scalar conservation law on a Riemannian manifold. In the case of the Torus, we show that, under a weak property of genuine non-linearity of the flux, the solution converges to its average value in $L^{p}$, $1\leq p<+\infty$. We give a partial result in the general case.
\end{abstract}

\medskip



\section{Introduction}\label{intro}


Let $M$ be a closed ({\it i.e.} compact, without boundary) Riemannian manifold of dimension $d$ with metric $g$. Let $\div$ be the divergence operator induced by $g$, defined in local coordinates by the formula $\div X=\partial_j X^j+\Gamma^j_{kj}X^k$, where $X$ is a vector field over $M$ and $\Gamma^i_{kj}$ the Cristoffel's symbols of $g$. Let $A(\cdot,u)$ be a one-parameter vector field over $M$ such that $A\in C^1(M\times\R,TM)$ and, for all $u\in\R$, $A(\cdot,u)$ is divergence free. We are interested in the long-time behavior of the entropy solution to the Cauchy Problem
\begin{eqnarray}
u_{t}+\div(A(x,u))=0, & t>0, x\in M,\label{scl}\\
u(x,0)= u_{0}(x), & x\in M.\label{ICscl}
\end{eqnarray}
Here the initial datum $u_0\in L^\infty(M)$. The Cauchy Problem~\refe{scl}-\refe{ICscl} has been studied by Amorim, Ben-Artzi, LeFloch in \cite{AmorimBenArtziLeFloch05,BenArtziLeFloch07}. More generally, first order scalar equation in non-divergence form on manifolds have been studied by Panov in \cite{Panov97}. The case of a manifold with boundary has also been addressed by Panov in \cite{Panov07}.
\medskip

Let $\vg$ be the measure on $M$ induced by the metric $g$. Without loss of generality, we suppose that the volume of $M$ is $1$.
For $x\in M,\xi\in\R$, we set $a(x,\xi)=(\partial_u A)(x,\xi)$. At fixed $\xi\in\R$, $a(\cdot,\xi)$ is a free-divergence field over $M$ whose flow $(\varphi^\xi_t)_{t\in\R}$ is a group of measure-preserving diffeomorphisms. We say that $a(\cdot,\xi)$ is ergodic if $(\varphi^\xi_t)$ is, that is to say the invariant sets of $(\varphi^\xi_t)$ are or full or zero measure. By the Ergodic Theorem, this is equivalent to the following statement: for all $v\in L^2(M)$,
\begin{equation*}
\ds\frac{1}{T}\ds\int_0^T v\circ\varphi^\xi_t dt\to \ds\int_M v d\vg 
\end{equation*}
in $L^2(M)$ when $t\to+\infty$. The notion of non-degeneracy that we introduce can be considered as a quantified hypothesis of ergodic character on the averages of $a$. To be more explicit, let 
\begin{equation*}
\ds L^2_0(M):=\left\{v\in L^2(M);\int_M v d\vg=0\right\}
\end{equation*}
be the space of square-integrable functions over $M$ with zero mean value.
We introduce the following definition of non-degeneracy of the flux.
\begin{definition}[Non-degeneracy condition] Let $E$ be a Borel subset of $\R$. We say that $A$ (or $a$) is non-degenerate on $E$ if the quantity
\begin{equation}
|\bar a|(T;E):=\ds\sup\left\{\ds\int_E \left|\ds\frac{1}{T}\int_0^T\<v\circ\varphi^\xi_{-t},v\circ\varphi^\zeta_{-t}\>_{L^2(M)}dt \right|d\xi\right\},
\label{Andeta}\end{equation}
where the supremum is taken over all $v\in L^2_0(M)$ with $\|v\|_{L^2(M)}=1$ and over $\zeta\in E$, tends to $0$ when $T\to+\infty$.
\label{def:And}\end{definition}

The notion of entropy solution to the Cauchy Problem~\refe{scl}-\refe{ICscl} is defined in \cite{AmorimBenArtziLeFloch05,BenArtziLeFloch07}. We will also give an equivalent definition in Definition~\ref{defeskin}. Our main result is the following one. 

\begin{theorem}[Long-time convergence of the entropy solution] {\ \\} 
Suppose that $M=\T^d$ is the $d$-dimensional torus and that $a$ is independent on $x$. Let $u_{0}\in L^{\infty}(M)$ and let $\bar u$ denote the mass (constant in time) of the entropy solution $u\in C(\R_{+};L^{1}(M))\cap L^{\infty}(\R_{+}\times M)$ of the Cauchy Problem~\refe{scl}-\refe{ICscl}, 
\begin{equation*}
\bar u=\int_M u(x,t) d\vg(x).
\end{equation*} 
If the flux $A$ is non-degenerate on a neighborhood of $\bar u$, then $u(t)\to \bar u$ in $L^p(M)$ for every $1\leq p<+\infty$ when $t\to+\infty$.
\label{thCVes}\end{theorem}

The result of convergence is given only in the case $M=\T^d$. The obstruction to the convergence in the general case is the lack of a result of compactness (see Section~\ref{sec:compactness}). However, a large part of the intermediate statements and results are more intuitive in the context of general Riemannian manifolds and we have kept this framework. Observe that, in the case of the torus, the non-degeneracy condition
\begin{equation*}
\forall \nu\in\S^{d}, |\{\xi\in E;(1,a(\xi))^T\cdot\nu\}|=0
\end{equation*}
implies the non-degeneracy condition given in Def.~\ref{def:And} (see Lemma~\ref{ndVSnd} at the end of the paper). Observe also that a condition of non-degeneracy is ne\-ces\-sary for the validity of the result. For example, if $M=\T^2$ is the $2$-dimensional torus, and $A(\xi)=\frac{1}{2}\xi^2 a_0$ where $a_0=(1, \alpha)^T\in\R^2$, then every function of the form $U(y-\alpha x)$ is a stationary solution to the  conservation law $u_t+\div A(u)=0$. 
\bigskip

In case of unbounded domains, the analysis of the long-time behaviour of the entropy solution includes in particular the study of the stability of shocks. We refer in particular to the work by Serre \cite{Serre04}. In the case of periodic domains, on which we focus, there are no travelling waves (under a non-degeneracy condition), but convergence to a constant as described in Theorem~\ref{thCVes}. As in previous results in the area, our proof consists in the study of the evolution over the $\omega$-limit set of a trajectory. This is the method used in \cite{DaFermos85,EEngquist93,ChenFrid99,ChenPerthame08}. The original result of Lax~\cite{Lax57} uses the Hopf-Lax formula. Let us also mention the references \cite{FeireislPetzeltova97}, \cite{Frid02} and \cite{AmadoriSerre06} for related work (respectively first-order scalar conservation law with memory, long-time behaviour in the almost periodic case and homogenization in periodic, forced scalar conservation law). We have learned after our work was finished that Chen and Perthame have written recently an article which have some similarity with ours (see \cite{ChenPerthame08}). By use of the kinetic formulation, they prove the long-time convergence to the mean value on the torus. Their result is valid for second order, possibly degenerate, scalar conservation law under a global non-degeneracy condition.  In comparison, our result is for first-order scalar conservation law but we use a non-degeneracy condition which is weaker: it is localized around the mean-value $\bar u$.
\bigskip

The proof of the convergence of the entropy solution $u$ stated in Theorem~\ref{thCVes} uses the kinetic formulation of the Cauchy Problem~\refe{scl}-\refe{ICscl}, which we give in Section~\ref{sec:kin}. In  Section~\ref{sec:freeT}, we prove a property of homogenization of the free transport equation $(\partial_t+a(x,\xi)\cdot\nabla)f=0$. In Section~\ref{sec:compactness} is given a result of compactness on the family $\{u(t);t>0\}$ (restricted to the case of the torus). At last in Section~\ref{sec:proof} we give the proof of Theorem~\ref{thCVes}.


\section{Kinetic formulation of conservation laws}\label{sec:kin}


The kinetic formulation of conservation laws dates back to the paper by Lions, Perthame, Tadmor~\cite{LionsPerthameTadmor94} and, retrospectively, can already be figured out in \cite{GigaMiyakawa83,Brenier83,PerthameTadmor91}. We refer to the book by Perthame~\cite{PerthameBook} concerning the subject. Although the theory has been addressed in the case $M=\R^N$ (or $M$ open subset of $\R^N$), the case of a closed Riemannian manifold is completely similar and therefore, instead of using the definition of entropy solution \`a la Kruzhkov and then give the equivalent kinetic formulation, we will at once define an entropy solution {\it via} the kinetic formulation. We leave it to the reader to verify that it coincides with the definition given in \cite{BenArtziLeFloch07}. We first introduce the definition of equilibrium functions.

For $\alpha,\xi\in\R$, define 
\begin{equation*}
\chi_\alpha(\xi)=\left\{\begin{array}{r l}
-1 &\mbox{if } \alpha<\xi<0,\\
1 &\mbox{if } 0<\xi<\alpha,\\
0 &\mbox{otherwise.}
\end{array}\right.
\end{equation*}
If $u\in L^1(M)$, $\chi_u$ is the {\it equilibrium function} associated to $u$. Notice that $\chi_u\in L^1(M\times\R)$, as shown by the following identities ($\alpha,\beta\in\R$):
\begin{equation}
\ds\int_\R(\chi_\alpha(\xi)-\chi_\beta(\xi))^+d\xi=(\alpha-\beta)^+,\quad
\ds\int_\R|\chi_\alpha(\xi)-\chi_\beta(\xi)|d\xi=|\alpha-\beta|.\label{intchi}
\end{equation}

\begin{definition}[Entropy solution] Let $u_0\in L^\infty(M)$. A function $u$ in $L^\infty(M\times(0,+\infty))$ is said to be an entropy solution to the Cauchy Problem~\refe{scl}-\refe{ICscl} if there exists a non-negative measure $m$ on $M\times[0,+\infty)\times\R$ such that, for all $\varphi\in C^1_c(M\times[0,+\infty)\times\R)$,
\begin{multline}
\ds\int_0^{\infty}\hskip -5pt\int_M\int_\R \chi_{u}(\partial_t+a(x,\xi)\cdot\nabla)\varphi d\xi d\vg(x) dt+\int_M\int_\R \chi_{u_0}\varphi(0) d\xi d\vg(x)\\
=\ds\int_0^{\infty}\hskip -5pt\int_M\int_\R \partial_\xi\varphi dm(t,x,\xi).
\label{weakeskin}\end{multline}
Additionally, $\pi_\# m$, the push-forward of $m$ by the projection $\pi\colon M\times[0,+\infty)\times\R\to\R$ on the $(x,t)$-coordinates satisfies $\pi_\# m\in L^1\cap L^\infty(\R)$, with
\begin{equation*}
\pi_\# m(\xi)\leq \|u_0\|_{L^1(M)},\quad \int_\R \pi_\# m(\xi) d\xi\leq\frac{1}{2}\|u_0\|_{L^2(M)}^2.
\end{equation*}
\label{defeskin}\end{definition}

Choose $\varphi(x,t,\xi)=\psi(x,t)\eta'(\xi)$, $\eta$ convex, $\psi\geq 0$. in \refe{weakeskin} and use the sign of the right hand-side:
\begin{equation*}
\ds\int_0^{\infty}\hskip -5pt\int_M\int_\R  \psi(x,t)\eta''(\xi) dm(t,x,\xi)\geq 0,
\end{equation*}  
to deduce the entropy formulation from the kinetic formulation. The converse process starts from the entropy formulation with (semi-) Kruzhkov's entropy to {\it define}
\begin{equation}
m=-[\partial_t(u-\xi)^+ +\div(\sgn_+(u-\xi)(A(x,u)-A(x,\xi)))]
\label{defmkin}\end{equation}
and check, by derivation with respect to $\xi$, that \refe{weakeskin} is satisfied.
\medskip

Ben-Artzi and LeFloch in \cite{BenArtziLeFloch07} show that there exists a unique entropy solution to the Cauchy Problem~\refe{scl}-\refe{ICscl} and that, besides, $u$ is continuous in time with values in $L^1(M)$. The continuity in time implies in particular that $\chi_u$ satisfies the weak equation with terminal time $T$:
\begin{multline}
\ds\int_0^{T}\hskip -5pt\int_M\int_\R \chi_{u}(\partial_t+a(x,\xi)\cdot\nabla)\varphi d\xi d\vg(x) dt+\int_M\int_\R \chi_{u_0}\varphi(0) d\xi d\vg(x)\\
-\int_M\int_\R \chi_{u(T)}\varphi(T) d\xi d\vg(x)=\ds\int_0^{T}\hskip -5pt\int_M\int_\R \partial_\xi\varphi dm(t,x,\xi),
\label{weakesT}\end{multline}
for every $\varphi\in C^1_c(M\times[0,T]\times\R)$.
\medskip

In \cite{BenArtziLeFloch07}, the author also prove the following contraction principle: if $u$ and $v$ are respectively the entropy solutions associated to the initial data $u_0$ and $v_0$ then 
\begin{equation}
\|(u(t')-v(t'))^+\|_{L^1(M)}\leq \|(u(t)-v(t))^+\|_{L^1(M)}\leq \|(u_0-v_0)^+\|_{L^1(M)}
\label{l1comparison}\end{equation}
for all $t'\geq t\geq 0$. This, in particular, implies the following maximum principle:
\begin{equation}
\|u(t)\|_{L^\infty(M)}\leq \|u_0\|_{L^\infty(M)}\mbox{ for all }t\geq 0.
\label{Pmax}\end{equation}

Note that, in Definition~\ref{defeskin}, we state that the measure $\pi_\# m$ is {\it finite} on $M\times[0,+\infty)$: here we emphasize the fact that the time interval is infinite; this estimate, uniform in time, is crucial in the analysis of the long-time behaviour of the solution. Such a property is obtained (at least formally) by integration with respect to $(x,t)$ in the definition~\refe{defmkin}: this gives $\pi_\# m(\xi)\leq \int_M (u_0-\xi)^+ d\vg(x)$.


\section{Long-time behavior in the free transport equation}\label{sec:freeT}


In this section we show that, under the property of non-degeneracy of $a$ defined in Def.~\ref{Andeta}, the density of the solution to the transport Problem 
\begin{equation}
\left\{\begin{array}{r l l l}
(\partial_{t}+a(x,\xi)\cdot\nabla_x)f(t,x,\xi)&=&0, & t>0, x\in M,\xi\in\R,\\
f(0,x,\xi)&=&f_0(x,\xi), & x\in M,\xi\in\R,
\end{array}\right.
\label{transportpur}\end{equation}
converges to its mean value. We first introduce some notations.

\begin{definition}[Density and mass] Let $f\in L^1(M\times\R)$ and let $E$ be a Borel subset of $\R$. Respectively the density and the mass of $f$ on $E$ are the quantities $u_E$ and $\bar u_E$ defined by
\begin{equation*}
u_E(x)=\ds\int_E f(x,\xi) d\xi,\quad x\in M,\quad \bar u_E=\ds\int_M u_E d\vg=\ds\int_{M\times E} f(x,\xi) d\xi d\vg.
\end{equation*}
When $E=\R$, these quantities are denoted $u$ and $\bar u$ and called respectively the density of $f$ and the mass of $f$.
\end{definition}

Since $M$ is of volume $1$ we also call $\bar u$ the mean-value of $u$. In what follows we denote by $\varphi^\xi_{t}$ the flow of $a(\cdot,\xi)$. The solution of the transport Problem~\refe{transportpur} will then be $f(t)\colon (x,\xi)\mapsto f_0(\varphi^\xi_{-t}x,\xi)$

\begin{proposition}[Homogenization in the free transport equation] Let $f_0\in L^2(M\times\R)$. Let $f\in C([0,+\infty);L^2(M\times\R))$ be the solution to the transport Problem \refe{transportpur}.
Let $E$ be a Borel subset of $\R$. If $A$ is non-degenerate on $E$, then the density of $f(t)$ on $E$ converges on average to the mass $\bar u_E$ when $t\to+\infty$:
\begin{equation}
\ds\frac{1}{T}\ds\int_{0}^T\|u_E(t)-\bar u_E\|^2_{L^2(M)} dt \leq |\bar a|(T;E)\|f_0\|^2_{L^2(M\times E)},
\label{estimT}\end{equation}
where $|\bar a|(T;E)$ is defined in \refe{Andeta}.
\label{propT}\end{proposition}

The decay described in \refe{estimT} above can be interpreted as a property of dispersion of the transport equation. We refer to te paper by Castella, Perthame~\cite{CastellaPerthame96} on the subject.
\medskip

\noindent {\bf Proof of Proposition~\ref{propT}:} by linearity of the transport equation, $u_E(x,t)-\bar u_E=\int_E \tilde f(t,x,\xi) d\xi$ where $\tilde f$ is the solution to the transport Problem~\refe{transportpur} with initial datum $\tilde f_0=f_0-\int_M f_0(x,\xi) d\vg$. Therefore, we can as well suppose that $f_0\in L^2(E;L^2_0(M))$ and that $\bar u_E=0$. By expanding the square in the $L^2$ norm, using Fubini's Theorem and the change of variable $x'=\varphi_{-t,\xi}x$, we have
\begin{equation*}
\|u_E(t)\|^2_{L^2(M)}= \int_{E\times E\times M} f_0(\varphi^\xi_{-t}(x),\xi)f_0(\varphi^{\zeta}_{-t}(x),\zeta) d\vg(x) d\zeta d\xi. 
\end{equation*}
Using Fubini's Theorem again, we obtain 
\begin{equation*}
\ds\frac{1}{T}\ds\int_{0}^T\|u_E(t)\|^2_{L^2(M)}dt= \int_{E}\int_E  \ds\frac{1}{T}\int_0^T\<f_0(\varphi^\xi_{-t}\cdot,\xi),f_0(\varphi^{\zeta}_{-t}\cdot,\zeta)\>_{L^2(M)} dt d\zeta d\xi.
\end{equation*}
The lemma then follows from the estimate
\begin{multline}
\ds\int_E\int_E \left|\ds\frac{1}{T}\int_0^T\<f(\varphi^\xi_{-t}\cdot,\xi),f(\varphi^\zeta_{-t}\cdot,\zeta)\>_{L^2(M)}dt \right|d\xi d\zeta\\
 \leq |\bar a|(T;E)\|f\|^2_{L^2(M\times E)},
\label{Andetaplus}\end{multline}
valid for all $f\in L^2(E;L^2_0(M))$. This inequality follows from \refe{Andeta} and elementary arguments. First, it is sufficient to prove it in case of separate variables. Indeed, the set of tensor functions $v\otimes k(x,\xi)=v(x)k(\xi)$, $v\in L^2(M),k\in L^2(E)$ is dense in $L^2(M\times E)$ and we notice that if $(v_n\otimes k_n)$ is a sequence of $L^2(M)\otimes L^2(E)$ converging in $L^2(M\times E)$ to a $f\in L^2(E;L^2_0(M))$, then we can suppose, upon subtracting its mean value to $v_n$, that $v_n\in L^2_0(M)$ for all $n$. Therefore, to prove \refe{Andetaplus}, we may as well suppose that $f=v\otimes k$, $v\in L^2_0(M)$, $k\in L^2(E)$. Then we have 
\begin{equation*}
\<f(\varphi^\xi_{-t}\cdot,\xi),f(\varphi^\zeta_{-t}\cdot,\zeta)\>_{L^2(M)}=k(\xi) k(\zeta) \<v\circ\varphi^\xi_{-t},v\circ\varphi^\zeta_{-t}\>_{L^2(M)}
\end{equation*}
for a.e. $\xi,\zeta\in E$. We use the elementary inequality $2|k(\xi)||k(\zeta)|\leq |k(\xi)|^2+|k(\zeta)|^2$ to obtain the bound
\begin{multline*}
\ds\int_E\int_E \left|\ds\frac{1}{T}\int_0^T\<f(\varphi^\xi_{-t}\cdot,\xi),f(\varphi^\zeta_{-t}\cdot,\zeta)\>_{L^2(M)}dt\right|d\xi d\zeta \leq\\
\ds\int_E\int_E |k(\zeta)|^2\left|\ds\frac{1}{T}\int_0^T \<v\circ\varphi^\xi_{-t},v\circ\varphi^\zeta_{-t}\>_{L^2(M)}\right| d\xi d\zeta
\end{multline*}
from which \refe{Andetaplus} follows by definition of $|\bar a|(T;E)$. This concludes the proof of Proposition~\ref{propT}. \qed


\section{Compactness of the orbit}\label{sec:compactness}


In this section, we state the following result of compactness.

\begin{theorem} Assume $M=\T^d$ is the $d$-dimensional torus and $a$ is independent on $x$. Let $u_0\in L^\infty(M)$ and let $u\in L^\infty(M\times(0,+\infty))\cap C([0,+\infty);L^1(M))$ be the entropy solution to the Cauchy Problem~\refe{scl}-\refe{ICscl}. Then the family $\{u(\cdot+t);t\geq 0\}$ is relatively compact in $L^1(M)$.
\label{th:compact}\end{theorem}

The result is classical in the theory of scalar conservation laws; we simply recall the arguments of the proof. We want to prove that the orbit is totally bounded, {\it i.e.} can be covered by a finite number of balls of arbitrary diameter in $L^1(M)$: first, by approximation and the property of $L^1$-contraction, we can suppose that $u_0\in BV(M)$. The invariance by translation of the problem and the property of $L^1$-contraction then show that
\begin{equation*}
h^{-1}\|u(t,\cdot+h)-u(t)\|_{L^1(M)}\leq h^{-1}\|u_0(\cdot+h)-u_0\|_{L^1(M)}\leq \|u_0\|_{BV(M)}
\end{equation*}
for any $h\in\R^d$, which yields the uniform bound $\|u(t)\|_{BV(M)}\leq\|u_0\|_{BV(M)}$, whence the result by compactness of the injection $L^1\cap BV(M)$ in $L^1(M)$.
\medskip


\section{Convergence of the entropy solution}\label{sec:proof}


In this section we prove the main result, Theorem~\ref{thCVes}. Let therefore $u_0\in L^\infty(M)$, let $\bar u$ be the mean-value over $M$ of $u_0$ and let $E$ be an open interval of $\R$ containing $\bar u$ such that $A$ is non-degenerate on $E$. Let $u$ be the entropy solution to \refe{scl}-\refe{ICscl}. By the maximum principle~\refe{Pmax}, we have
$-R\leq u\leq R$ a.e. where $R:=\|u_0\|_{L^\infty(M)}$. If $\lambda\in\R$, then $u+\lambda$ is also solution to \refe{scl}-\refe{ICscl} where $u_0$ has been replaced by $u_0+\lambda$ and $A(u,x)$ by $A(u-\lambda,x)$. These modifications do not affect the hypothesis of non-degeneracy of $A$, therefore we will suppose, without loss of generality, that $u$ and $u_0$ are non-negative: $0\leq u\leq R$ and $0\leq u_0\leq R$ a.e. We will also discard the trivial cases $\bar u=0$ or $\bar u=R$. Eventually, we will suppose in a first step that $E=(0,R)$ ($A$ is non-degenerate everywhere).
\medskip

{\bf Step 1. $A$ non-degenerate everywhere.} Let $(t_n)$ be an increasing sequence of time steps which tends to $+\infty$. Let $\eps>0$. Fix $T_1>0$ such that $|\bar a|(T_1,E)\|u_0\|_{L^1(M)}<\eps^2$. 
Up to a subsequence, we can suppose that all the intervals $[t_n,t_n+T_1]$, $n\in\N$, are mutually disjoint.
Since $m$ is finite over $M\times[0,+\infty)\times[0,R]$, we have then
\begin{equation}
\ds\lim_{n\to+\infty} m(M\times[t_n,t_n+T_1]\times[0,R])=0.
\label{killm}\end{equation}
Define $u_n\in C([0,T_1];L^1(M))$ by $u_n(x,t)=u(x,t+t_n)$. By the compactness result of Theorem~\ref{th:compact}, there exists a subsequence still denoted $(t_n)$ and $u_{0,\infty}\in L^\infty(M)$ such that $u_n(0)\to u_{0,\infty}$ in $L^1(M)$. Since $f_n:=\chi_{u_n}\in [0,1]$ a.e. on $M\times[0,T_1]\times\R$, there exists $f_\infty\in L^\infty(M\times[0,T_1]\times\R)$, $f_\infty\in [0,1]$ a.e., such that $f_n\rightharpoonup f_\infty$ in $L^\infty(M\times[0,T_1]\times\R)$-weak-*. 
Additionnally, $\partial_\xi f_n=\delta_0(\xi)-\nu^n_{t,x}(\xi)$ where $\nu^n_{t,x}$ is the Young measure with support in $[0,R]$ defined by $\nu^n_{t,x}=\delta_{u_n(t,x)}$. Consequently, up to a new subsequence, we can suppose that $\partial_\xi f_\infty=\delta_0(\xi)-\nu_{t,x}$ where $\nu_{t,x}$ is a Young measure supported in $[0,R]$. The equation satisfied by $f_\infty$ is the following one: for $\varphi\in C^1_c(M\times[0,T_1)\times\R)$, we have, introducing $\varphi_n(t,x,\xi)=\varphi(t-t_n,x,\xi)$, $f_n=\chi_{u_n}$, 
\begin{equation*}
\ds\int_0^{T_1}\hskip -5pt\int_M\int_0^R f_n(\partial_t+a\cdot\nabla)\varphi=\ds\int_{t_n}^{t_n+T_1}\hskip -5pt\int_M\int_0^R f(\partial_t+a\cdot\nabla)\varphi_n.
\end{equation*}
By \refe{weakesT}, we compute explicitely the right hand-side and obtain
\begin{equation}
\ds\int_0^{T_1}\hskip -5pt\int_M\int_0^R f_n(\partial_t+a\cdot\nabla)\varphi+\int_M\int_0^R \chi_{u_n(0)}\varphi(0)=\ds\int_{t_n}^{t_n+T_1}\hskip -7pt\int_M\int_0^R m\partial_\xi\varphi_n.
\label{proof0}\end{equation}

The bound $\|\partial_\xi \varphi_n\|_{L^\infty([t_n,t_n+T_1]\times M\times [0,R])}\leq\|\partial_\xi \varphi\|_{L^\infty([0,T_1]\times M\times [0,R])}$ and \refe{killm} show that the right-hand side of \refe{proof0} converges to $0$ when $n\to+\infty$:
\begin{equation*}
\ds\int_0^{T_1}\hskip -5pt\int_M\int_0^R f_\infty(\partial_t+a\cdot\nabla)\varphi+\int_M\int_0^R \chi_{u_{\infty,0}}\varphi(0)=0.
\end{equation*}
Consequently, $f_\infty$ is solution to the free transport equation $(\partial_t+a\cdot\nabla)f=0$ with initial datum $\chi_{u_{\infty,0}}$. Using test-functions of the form $\varphi(t,x)\eta'(\xi)$ for regular convex function $\eta$, we obtain that $\nu_{t,x}$ is a measure-valued (is-)entropy solution to the scalar conservation law~\refe{scl} with initial datum $\delta_{u_{0,\infty}}$. By Theorem~5.3 in \cite{BenArtziLeFloch07}, which asserts the uniqueness of measure-valued entropy solution with Dirac mass initial datum, we have $\nu_{x,t}=\delta_{u_\infty(x,t)}$ where $u_\infty$ is the entropy solution to \refe{scl} with initial datum $u_{0,\infty}$. Coming back at the kinetic level, we obtain that $f_\infty$ is an equilibrium function: $f_\infty=\chi_{u_\infty}$. This also implies that $(u_n)$ is converging to $u_\infty$ in $L^1(M\times[0,T_1])$ strongly.
\medskip

Since $f_\infty$ has the same mass as $f_n(t)$, $\bar u_\infty=\bar u$, Proposition~\ref{propT} with, we recall, $E=(0,R)$ here, gives
\begin{equation*}
\ds\frac{1}{T_1}\int_0^{T_1}\|u_\infty(t)-\bar u\|_{L^2(M)}^2 dt\leq|\bar a|(T_1,E)\|\chi_{u_{\infty,0}}\|_{L^2(M\times[0,R])}^2.
\end{equation*}
Since $\|\chi_{u_{\infty,0}}\|_{L^2(M\times[0,R])}^2=\|u_{\infty,0}\|_{L^1(M)}\leq \|u_0\|_{L^1(M)}$, we obtain
\begin{equation}
\ds\frac{1}{T_1}\int_0^{T_1}\|u_\infty(t)-\bar u\|_{L^2(M)}^2 dt<\eps^2.
\label{proof2}\end{equation}
Then, Jensen's Inequality and the inequality $\|v\|_{L^1(M)}\leq \|v\|_{L^2(M)}$ show finally that
\begin{equation}
\ds\frac{1}{T_1}\int_0^{T_1}\|u_\infty(t)-\bar u\|_{L^1(M)} dt<\eps.
\label{proof22}\end{equation}
By the $L^1$-contraction/comparison principle~\refe{l1comparison}, the map $t\mapsto\|u(t)-\bar u\|_{L^1(M)}$ is non-increasing. Therefore it has a limit when $t\to+\infty$ and, furthermore, it satisfies for every $n\in\N$,
\begin{equation*}
\|u((T_1+t_n))-\bar u\|_{L^1(M)}\leq\ds\frac{1}{T_1}\ds\int_{0}^{T_1}\|u_n(t)-\bar u\|_{L^1(M)}dt.
\end{equation*}
At the limit $[n\to+\infty]$, \refe{proof22} gives
\begin{equation*}
\ds\lim_{t\to+\infty}\|u(t)-\bar u\|_{L^1(M)}\leq\frac{1}{T_1}\ds\int_{0}^{T_1}\|u_\infty(t)-\bar u\|_{L^1(M)}dt<\eps.
\end{equation*}
Since $\eps$ is arbitrary, $\lim_{t\to+\infty}\|u(t)-\bar u\|_{L^1(M)}=0$.
\medskip

{\bf Step 2. $A$ non-degenerate in the neighborhood of $\bar u$.} Let us now turn to the general case, where $E$ is a neighborhood of $\bar u$ and not necessarily the whole interval $(0,R)$. Without loss of generality, we can suppose $E=[\bar u,\zeta]$, $\zeta\in(\bar u,R)$ (recall that we have supposed $0<\bar u<R$). Fix $\eps_0>0$ (possibly depending on $E$ and $\bar u$), fix $\eps<\min(\eps_0,1)$ and $T_1>0$ such that $|\bar a|(T_1,E)\|u_0\|_{L^1(M)}<\eps^4$. The same reasoning as above gives, instead of \refe{proof2}, the bound
\begin{equation*}
\ds\frac{1}{T_1}\int_0^{T_1}\|u_{\infty,E}(t)-\bar u_{\infty,E}\|_{L^2(M)}^2 dt<\eps^4,
\end{equation*}
where $u_{\infty,E}$ and $\bar u_{\infty,E}$ are respectively the density and mass of $f_\infty$ in $E$.
By Jensen's Inequality, and the inequality $\|v\|_{L^1(M)}\leq \|v\|_{L^2(M)}$, we obtain the estimate
\begin{equation}
\ds\frac{1}{T_1}\int_0^{T_1}\|u_{\infty,E}(t)-\bar u_{\infty,E}\|_{L^1(M)} dt\leq\eps^2.
\label{qBGK3}\end{equation}
To conclude in the same way as above ({\it i.e.} by monotony of $t\mapsto\|u(t)-\bar u\|_{L^1(M)}$), it is therefore sufficient to show that the norm $\|u_{\infty,E}(t)-\bar u_{\infty,E}\|_{L^1(M)}$ gives a control on the norm $\|u_\infty(t)-\bar u_\infty\|_{L^1(M)}$. This is the content of the following lemma that we apply to $v=u_\infty(t)$ and $f=f_\infty(t)$ (notice that the statement and result of the lemma are time-independent).


\begin{lemma} Let $v\in L^\infty(M)$, $0\leq v\leq R$ a.e. and let $f=\chi_v$ be the associated equilibrium function. Let $E$ be an interval of the form $E=[\bar v,\zeta]$, $\zeta\in(\bar v,R)$. Then there exists $\delta_0=\delta_0(|E|,\bar v)>0$ such that, for every $0<\delta<\delta_0$,
\begin{equation}
\|v-\bar v\|_{L^1(M)}\leq \delta+8\left(1+R\delta^{-1}\right)\|v_E-\bar v_E\|_{L^1(M)}.
\label{eqndloc}\end{equation}
\label{ndloc}\end{lemma}

Admit the lemma for the moment. Since $\bar u_\infty=\bar u$ and thus $\delta_0(|E|,\bar u_\infty)=\delta_0(|E|,\bar u)$ are independent on time, we can sum over $t\in[0,T_1]$ the inequality \refe{eqndloc} where $v=u_\infty(t)$; by \refe{qBGK3} we obtain the estimate
\begin{equation*}
\ds\frac{1}{T}\ds\int_0^T\|u_\infty(t)-\bar u\|_{L^1(M)}dt< \delta+8\left(1+R\delta^{-1}\right)\eps^2.
\end{equation*} 
Having chosen retrospectively $\eps_0=\delta_0$, we see we can take $\delta=\eps$ in the inequality above to get the following bound, similar to \refe{proof22} in Step~1:
\begin{equation*}
\ds\frac{1}{T}\ds\int_0^T\|u_\infty(t)-\bar u\|_{L^1(M)}dt< (9+8R)\eps.
\end{equation*}
We then derive, as in Step~1, the estimate $\ds\lim_{t\to+\infty}\|u(t)-\bar u\|_{L^1(M)}<(9+8R)\eps$ and conclude that $\lim_{t\to+\infty}\|u(t)-\bar u\|_{L^1(M)}=0$. \qed
\medskip


{\bf Proof of Lemma~\ref{ndloc}:} Note that since $E=[\bar v,\zeta]$
$$
v_E(x)=\min\{\zeta,v(x)\}-\bar v.
$$
We deduce the following formula
\begin{equation*}
(v(x)-\bar v)^+=v_E(x)+(v(x)-\zeta)^+
\end{equation*}
and will use more specifically the three following implications:
\begin{itemize}
\item If $v_E(x)>0$, then $v(x)\geq v_E(x)+\bar v$,
\item if $0\leq v_E(x)<|E|$, then $(v(x)-\bar v)^+\leq v_E(x)$,
\item If $0<v_E(x)<|E|$, then $|v(x)-\bar v-\bar v_E(x)|=|v_E(x)-\bar v_E(x)|$.
\end{itemize}
\medskip
 
We first show that $\bar v_E$ stays strictly lower than $|E|$: we have $v\geq v_E+\bar v$ whenever $v_E\geq\lambda>0$, and therefore
\begin{equation*}
\bar v\geq\int_{v_E>\lambda}v\geq\int_{v_E>\lambda}(v_E+\bar v).
\end{equation*}
Integrating over $\lambda\in[0,|E|]$, we get $\bar v|E|\geq \int_M v_E^2+\bar v v_E$ and, by Jensen's Inequality, $\bar v|E|\geq \bar v_E^2+\bar v v_E$.
Consequently, $\bar v_E$ is smaller than the root $\kappa$ in $[0,|E|]$ of the equation $\kappa^2=(|E|-\kappa)\bar v$. Since this root cannot be $|E|$, it is strictly smaller than $|E|$: $\kappa\leq |E|-c(|E|,\bar v)$ where $c(|E|,\bar v)>0$. We obtain:
\begin{equation}
\bar v_E\leq |E|-c(|E|,\bar v).
\label{notallE}\end{equation}
Set $\delta_0=\min(|E|/4,c(|E|,\bar v))$ and fix $\delta\in (0,\delta_0)$. If $0<v_E(x)\leq 2\delta$, then $0\leq v_E(x)<|E|$ and thus $(v-\bar v)^+(x)\leq v_E(x)\leq 2\delta$. This shows that
\begin{equation}
\ds\int_M (v-\bar v)^+\leq 2\delta+\ds\int_{\{v_E> 2\delta\}} (v-\bar v)^+\leq 2\delta+R|\{v_E> 2\delta\}|.
\label{K10}\end{equation}
If $\bar v_E\leq\delta$, then, by Chebychev's Inequality,
\begin{equation*}
|\{v_E> 2\delta\}|\leq |\{v_E-\bar v_E> \delta\}|\leq \delta^{-1}\|v_E-\bar v_E\|_{L^1(M)}.
\end{equation*}
Since $\|v-\bar v\|_{L^1(M)}=2\|(v-\bar v)^+\|_{L^1(M)}$, we get the estimate
\begin{equation}
\|v-\bar v\|_{L^1(M)}\leq 4\delta+2R\delta^{-1}\|v_E-\bar v_E\|_{L^1(M)},
\label{K1}\end{equation}
provided $\bar v_E\leq\delta$. In the case $\bar v_E>\delta$ we proceed slightly differently: we give an estimate on the quantity
\begin{equation*}
\int_M |v-\bar v-\bar v_E|
\end{equation*} 
via the partition $M=M_1\cup M_2$, $M_1=\{|v_E-\bar v_E|<\delta\}$, $M_2=\{|v_E-\bar v_E|>\delta\}$. On $M_1$, we have $0<v_E<|E|$ by \refe{notallE} and the choice of $\delta_0$, hence
$|v-\bar v-\bar v_E|=|v_E-\bar v_E|$. This shows that
\begin{equation*}
\int_{M_1} |v-\bar v-\bar v_E|\leq \|v_E-\bar v_E\|_{L^1(M)}.
\end{equation*}
On the other hand, the sum over $M_2$ can be bounded as above in \refe{K1}:
\begin{equation*}
\int_{M_2} |v-\bar v-\bar v_E|\leq R\delta^{-1}\|v_E-\bar v_E\|_{L^1(M)}.
\end{equation*}
We obtain therefore $\|v-\bar v-\bar v_E\|_{L^1(M)}\leq\left(1+R\delta^{-1}\right)\|v_E-\bar v_E\|_{L^1(M)}$.
Since $\bar v$ is the mean value of $v$, we also have 
\begin{equation*}
0\leq \bar v_E=-\int_M v-\bar v-\bar v_E\leq \|v-\bar v-\bar v_E\|_{L^1(M)}.
\end{equation*}
This shows that $\|v-\bar v\|_{L^1(M)}\leq 2\left(1+R\delta^{-1}\right)\|v_E-\bar v_E\|_{L^1(M)}$ when $\bar v_E>\delta$. By \refe{K1}, we conclude that, independently on the sign of $\bar v_E-\delta$, we have $\|v-\bar v\|_{L^1(M)}\leq 4\delta+2\left(1+R\delta^{-1}\right)\|v_E-\bar v_E\|_{L^1(M)}$. Replace $\delta$ by $\delta/4$ to obtain \refe{eqndloc}. \qed
\bigskip


Notice that, apart from the argument of relative compactness of the orbit, the proof of Theorem~\ref{thCVes} is not restricted to the case of the torus. We can therefore state the following result.

\begin{proposition} Let $u_{0}\in L^{\infty}(M)$ and let $\bar u$ denote the mass (constant in time) of the entropy solution $u\in C(\R_{+};L^{1}(M))\cap L^{\infty}(\R_{+}\times M)$ of the Cauchy Problem~\refe{scl}-\refe{ICscl}, 
\begin{equation*}
\bar u=\int_M u(x,t) d\vg(x).
\end{equation*} 
Assume that the orbit $\{u(t);t\geq 0\}$ is relatively compact in $L^1(M)$. If the flux $A$ is non-degenerate on a neighborhood of $\bar u$, then $u(t)\to \bar u$ in $L^p(M)$ for every $1\leq p<+\infty$ when $t\to+\infty$.
\label{propCVes}\end{proposition}

A possible way to prove the compactness of the orbit would be the extension of the averaging lemmas to the case of Riemannian manifold, under the non-degeneracy hypothesis of Def.~\ref{def:And}. The proof of the averaging lemmas being essentially based on the Fourier Transform (case of $\R^d$) or Fourier Series (case of $\T^d$), this seems out of reach for the moment.
\bigskip

To complete Theorem~\ref{thCVes}, we prove the following lemma.

\begin{lemma} Let $E$ be a Borel subset of $\R$. Assume the non-degeneracy condition
\begin{equation}
\forall \nu\in\S^{d}, |\{\xi\in E;(1,a(\xi))^T\cdot\nu\}|=0.
\label{Andclassic2}\end{equation}
Then $A$ is non-degenerate on $E$ in the sense of Definition~\ref{def:And}.
\label{ndVSnd}\end{lemma}

{\bf Proof:} Denote by $(e_n)_{n\in\Z^d}$ the Fourier orthonormal basis on $\T^d$, $e_n(x)=\exp(2i\pi n\cdot x)$. Since $\varphi^\xi_{-t}(x)=x+ta(\xi)$ is a translation, we have 
\begin{equation*}
\<e_n\circ\varphi^\xi_{-t},e_m\circ\varphi^\xi_{-t}\>_{L^2(M)}=\exp(2i\pi t(m\cdot a(\zeta)-n\cdot a(\xi)))\delta_{n,m}
\end{equation*}
for $n,m\in\Z^d$. This shows that $|\bar a|(T;E)=\sup_{n\in\Z^d\setminus\{0\}}|\bar a|_n(T;E)$ for
\begin{equation*}
|\bar a|_n(T;E):=\ds\sup\left\{\ds\int_E \left|\ds\frac{1}{T}\int_0^T\<e_n\circ\varphi^\xi_{-t},e_n\circ\varphi^\zeta_{-t}\>_{L^2(M)}dt \right|d\xi\right\},
\end{equation*}
where the supremum is taken over $\zeta\in E$. We compute
\begin{equation}
\ds\left|\frac{1}{T}\int_0^T\<e_n\circ\varphi^\xi_{-t},e_n\circ\varphi^\zeta_{-t}\>_{L^2(M)}dt\right|
=\left|\frac{\sin(\pi T n\cdot[a(\xi)-a(\zeta)])}{\pi T n\cdot[a(\xi)-a(\zeta)]}\right|
\label{explambda}\end{equation}
in case $n\cdot[a(\xi)-a(\zeta)]\not=0$ (value $1$ otherwise). For $1\geq \gamma>0$, set 
\begin{equation*}
\eps(\gamma):=\ds\sup_{\alpha\in\R,\beta\in\S^{d-1}}|\{\xi\in E;|\beta\cdot a(\xi)-\alpha|\leq\gamma\}|.
\end{equation*}
Notice that the supremum in $\eps(\gamma)$ can be taken over $K=[-M-1,M+1]\times\S^{d-1}$ where $M:=\|a\|_{L^\infty(E)}$.
By \refe{Andclassic2}, $\left|\left\{\xi\in E;\beta\cdot a(\xi)-\alpha=0\right\}\right|=0$ for all $\beta\in\S^{d-1},\forall\alpha\in\R$. In particular the function
$$
(\alpha,\beta,\gamma)\mapsto \left|\left\{\xi\in E;|\beta\cdot a(\xi)-\alpha|\leq\gamma\right\}\right|
$$
is continuous with respect to $\gamma$ and, at fixed $\alpha,\beta$, converges to $0$ when $\gamma\to 0$. Besides, the convergence is monotone, hence uniform with respect to $(\alpha,\beta)$ in the compact $K$. We conclude that $\eps(\gamma)\to 0$ when $\gamma\to 0$. 

Let $\eps>0$. There exists $\gamma>0$ such that $\eps(\gamma)<\eps$. Let $n\in\Z^d\setminus\{0\}$ and $\zeta\in E$. Set $E(n,\zeta)=\left\{\xi\in E;|n(\cdot a(\xi)-a(\zeta))|\leq\gamma\right\}$ and consider the sums
\begin{equation*}
\int_F \left|\frac{\sin(\pi T n\cdot[a(\xi)-a(\zeta)])}{\pi T n\cdot[a(\xi)-a(\zeta)]}\right| d\xi
\end{equation*}
for $F=E(n,\zeta)$ and $F=E\setminus E(n,\zeta)$. When $F=E(n,\zeta)$, the argument is bounded (by $1$) and the domain of integration is of measure less than $\eps$, thus the sum is less than $\eps$. Wen $F=E\setminus E(n,\zeta)$, the argument is bounded by $|\pi T \gamma|^{-1}$, thus the sum is less than $\frac{|E|}{\pi T \gamma}$. We conclude that 
\begin{equation*}
|\bar a|(T;E)<\eps+\frac{|E|}{\pi T \gamma},
\end{equation*} 
hence $|\bar a|(T;E)<2\eps$ for $T$ large enough. \qed

\def\cprime{$'$} \def\cprime{$'$}
\providecommand{\bysame}{\leavevmode\hbox to3em{\hrulefill}\thinspace}
\providecommand{\MR}{\relax\ifhmode\unskip\space\fi MR }
\providecommand{\MRhref}[2]{%
  \href{http://www.ams.org/mathscinet-getitem?mr=#1}{#2}
}
\providecommand{\href}[2]{#2}

\end{document}